\documentclass [11pt, a4paper, twoside, reqno, final]{amsart}

\usepackage[PS]{diagrams}
\usepackage{enumerate}
\usepackage{hyperref}

\theoremstyle{plain}
\newtheorem * {theorem}{Theorem}
\newtheorem * {lemma}{Lemma}
\newcommand{\nbd}{\nobreakdash}
\newcommand{\bR}{{\mathbb R}}
\newcommand{\bZ}{{\mathbb Z}}
\newcommand{\bN}{{\mathbb N}}
\newcommand{\iso}{\cong}
\newcommand{\tensor}{\otimes}
\newcommand{\id}{\mathrm{id}}

\newcommand{\pp}{\mathfrak{P}}
\newcommand{\ppgr}{\mathfrak{Pgr}}
\newcommand{\gr}{\mathrm{gr}}

\begin{document}

\title[$K$\nbd-theory of graded rings]{A note on the graded
  $K$\nbd-theory\\of certain graded rings}

\author{Thomas H\"uttemann}

\address{Thomas H\"uttemann, Queen's University Belfast, School of
  Mathematics and Physics, Pure Mathematics Research Centre, Belfast
  BT7~1NN, Northern Ireland, UK}

\email{t.huettemann@qub.ac.uk}

\urladdr{http://huettemann.zzl.org/}



\thanks{This work was supported by the Engineering and Physical
  Sciences Research Council [grant number EP/H018743/1].}

\begin{abstract}
  Following ideas of \textsc{Quillen} we prove that the graded
  $K$-theory of a $\bZ^n$\nbd-graded ring with support contained
  in a pointed cone is entirely determined by the $K$\nbd-theory of
  the subring of degree~$0$ elements.
\end{abstract}

\maketitle

Modules will be right modules, and the adjective ``graded'' refers to
``$\bZ^n$-graded'' throughout. Suppose that $R = \bigoplus_{a \in
  \bZ^n} R_a$ is a graded ring such that there exists a pointed
polyhedral cone $C \subset \bR^n$ containing the support of~$R$. That
is, $C$~is of the form $C = \mathrm{cone} \{v_1,\, v_2,\, \cdots,\,
v_\ell\}$ with $v_i \in \bR^n$, and there exists a non-zero linear
form~$\gamma$ on~$\bR^n$ with $\gamma (C) \subseteq [0,\infty[$ and $C
\cap \ker \gamma = \{0\}$. This implies that the hyperplanes
$\gamma^{-1} (\alpha)$ intersect~$C$ in a bounded set if $\alpha \geq
0$, and do not meet~$C$ for $\alpha < 0$. We will assume $\dim C = n$
(equivalently, $\mathrm{int}\, C \ne \emptyset$) as can be done
without loss of generality.

By slightly varying the linear form we can achieve that the
restriction $\gamma|_{\bZ^n}$ is injective; in fact, the set of linear
forms such that the $n$~real numbers $\gamma (e_1),\, \gamma(e_2),\,
\cdots,\, \gamma(e_n)$ are $\mathbb{Q}$\nbd-linearly independent
($e_i$~the $i$th unit vector) is dense in the set of all linear forms.
We will assume that we have chosen such a form~$\gamma$ from now
on. The convention
\[a \leq b \ :\Leftrightarrow \ \gamma(a) \leq \gamma(b)\] defines a
total order on~$\bZ^n$ such that for any $v \in \bR^n$ the set $(v +
C) \cap \bZ^n$ is well-ordered.

We define the ideal $R_+ = \bigoplus_{a > 0} R_a$, and identify~$R_0$
with~$R/R_+$.  A graded $R$\nbd-module is called {\it bounded below\/}
if its support is bounded below with respect to the above order. It is
easy to see that {\it if $M$~is a bounded below graded $R$\nbd-module
  with $M / MR_+ = 0$, then $M = 0$}. Any finitely generated graded
$R$\nbd-module is bounded below as its support is contained in a
translate of~$C$.

For a graded ring~$S$ let $\ppgr(S)$ denote the category of finitely
generated graded projective $S$\nbd-modules and their degree
preserving homomorphisms. Considering $R_0$ as a graded ring
concentrated in degree~$0$, we define a functor
\[T \colon \ppgr(R) \rTo \ppgr(R_0) \ , \quad P \mapsto P \tensor_R
R_0 = P / PR_+\ .\] With the remarks from the previous paragraph, we
can show just as in \cite[p.~637]{Bass-K}:

\begin{lemma}
  For a finitely generated graded projective $R$\nbd-module~$P$ there is
  a non-canonical isomorphism of graded $R$\nbd-modules
  \[P \iso T(P) \tensor_{R_0} R = \bigoplus_{b \in \bZ^n} T(P)_b
  \tensor_{R_0} R(-b)\] where the grading is such that
  \[P_a \iso \bigoplus_{b \in \bZ^n} T(P)_b \tensor_{R_0} R(-b)_a \
  .\]
\end{lemma}
More precisely, the short exact sequence of graded $R$\nbd-modules
\begin{equation}
  \label{eq:ses}
  0 \rTo \ker \pi \rTo P \rTo^\pi P / PR_+ = T(P) \rTo 0
\end{equation}
(with $\pi$~the canonical projection), when considered in the category
of graded $R_0$\nbd-modules, admits a degree-preserving splitting
$\sigma \colon T(P) \rTo P$ such that
\[T(P) \tensor_{R_0} R \rTo P \ , \quad x \tensor r \mapsto \sigma(x)
\cdot r\] is an isomorphism of graded $R$\nbd-modules. The splitting
map exists as $T(P)$~is projective as an $R_0$\nbd-module.

Write $\pp(R_0)$ for the category of finitely generated projective
$R_0$\nbd-modules, and consider $\pp(R_0)$ and $\ppgr(R)$ as exact
categories with the usual notion of exact sequences (exact in the
category of $R_0$\nbd-modules). We use the notation
\[K_i(R_0) = K_i(\pp(R_0)) \qquad \text{and} \qquad K_i^\gr (R) =
K_i(\ppgr(R))\] for the \textsc{Quillen} $K$\nbd-groups, and
abbreviate $\bZ[t_1,\, t_1^{-1},\, t_2,\, t_2^{-1},\, \cdots,\, t_n,\,
t_n^{-1}]$ by~$L$. The graded $K$\nbd-group $K_i^\gr(R)$ is a right
$L$\nbd-module with the action of~$t_i$ induced by the translation
functor $\ppgr(R) \rTo \ppgr(R),\ M \mapsto M(e_i)$ where $e_i$~is the
$i$th unit vector, and $M(e_i)_a = M_{e_i+a}$.

\begin{theorem}
  There is an isomorphism of $L$\nbd-modules
  \[K_i (R_0) \tensor L \rTo^\iso K_i^\gr (R) \ , \quad x \tensor 1
  \mapsto (\,\cdot\, \tensor_{R_0} R)_* (x) \ .\]
\end{theorem}

The Theorem and its proof are variations of \cite[Proposition,
p.~107--108]{Quillen-K}. For $a \in \bZ^n$ and $P \in \ppgr(R)$ let
$F^aP$ denote the graded submodule generated by the elements of~$P_d$
for $d \leq a$. Choose a non-zero integral vector~$v$ in the interior
of~$C$, and write $\ppgr_k(R)$ for the full subcategory of~$\ppgr(R)$
of modules~$P$ with support contained in~$-kv + C$ such that $F^{-kv}P
= 0$ and $F^{kv}P = P$.

Fix $P \in \ppgr_k(R)$ and $a \in \bZ^n$. By the Lemma, there is a
non-canonical isomorphism of $R_0$\nbd-modules
\[P_a \iso \bigoplus_b T(P)_b \tensor_{R_0} R(-b)_a \ .\] As $P$~is
finitely generated almost all summands are trivial. If $b \notin
-kv+C$ then $b$~is not in the support of~$P$ so that $T(P)_b$ and
hence the corresponding $b$\nbd-summand are trivial. If $b > a$ then
$R(-b)_a = 0$ so that the corresponding $b$\nbd-summand is trivial
again.

We now restrict attention to the $b$\nbd-summand for fixed $b \in
\bZ^n$. For~$P_c$ to be non-zero we need that $R(-b)_c = R_{c-b} \ne
0$ which necessitates $c-b \in C$ or, equivalently, $c \in b+C$. In
other words, for any $a \in \bZ^n$ the $R$\nbd-module $F^a \big(T(P)_b
\tensor R(-b)\big)$ is generated by the elements of $T(P)_b \tensor
R(-b)_{c}$ for $c \in (b+ C) \cap \bZ^n$ with $c \leq a$. Such~$c$
exist if and only if $a \geq b$, and we have $F^a \big(T(P)_b \tensor
R(-b)\big) = T(P)_b \tensor R(-b)$ in this case. If $a < b$ then $F^a
\big(T(P)_b \tensor R(-b)\big) = 0$.

As the assignment $F^a \colon P \mapsto F^aP$ is compatible with the
direct sum decomposition, this implies that there is a non-canonical
isomorphism
\begin{equation}
  \label{eq:noncan-Fa}
  F^a P \iso \bigoplus_{\begin{subarray}{l}b \in -kv+C \\ b \leq
      a\end{subarray}} T(P)_b \tensor_{R_0} R(-b)
\end{equation}
for any~$a$. It follows that $F^aP$ is finitely generated projective,
and that $F^a$~is an exact endo-functor of~$\ppgr_k(R)$.

List in ascending order all integral points
\[-kv = a_0 < a_1 < a_2 \ldots < a_m = kv\] in~$-kv+C$ with $a_i \leq
kv$. This gives rise to an admissible filtration
\begin{equation}
  \label{eq:filtration}
  0 = F^{a_0} \subseteq F^{a_1} \subseteq \ldots \subseteq F^{a_m} =
  \id_{\ppgr_k(R)}
\end{equation}
of exact endo-functors of~$\ppgr_k(R)$. Now for $a = a_j$ the
$b$\nbd-summand in~\eqref{eq:noncan-Fa} is trivial unless $b = a_i$
for some $i \leq j$. Hence the $b$\nbd-summands of $F^{a_j}P$
and~$F^{a_{j-1}}P$ agree except when~$b = a_j$ where the latter is
trivial. It follows that $F^{a_j}P/F^{a_{j-1}}P \iso T(P)_{a_j}
\tensor R(-a_j)$.

Surprisingly this last isomorphism is canonical. Indeed, the
isomorphism from the Lemma depends on a choice of splitting~$\sigma$
of the sequence~\eqref{eq:ses}, mapping an element $x \in T(P)_{a_j}$
to $\sigma(x) \in P_{a_j}$. Given another choice~$\sigma^\prime$ of a
splitting we observe that $(\sigma - \sigma^\prime)(x) \in \ker(\pi) =
PR_+$ so that $(\sigma - \sigma^\prime)(x) \in F^{a_{j-1}} P$. In
other words, passing to the quotient eliminates precisely the
indeterminacy of the isomorphism.

It follows from the previous paragraph that there are canonical
isomorphisms of functors
\begin{equation}
  \label{eq:iso}
  F^{a_j} P / F^{a_{j-1}} P \iso T(P)_{a_j} \tensor_{R_0} R(-a_j)
\end{equation}
for $1 \leq j \leq m$. The quotient functors $F^{a_j} / F^{a_{j-1}}$
are thus seen to be exact, and we can define the homomorphism
\[\Phi \colon \bigoplus_{j=0}^m K_i (R_0) \tensor t^{a_j} \rTo K_i
\big(\ppgr_k (R)\big) \ , \quad x \tensor t^{a_j} \mapsto G^j_*x\]
where $G^j \colon \pp(R_0) \rTo \ppgr_k(R),\ P \mapsto P \tensor_{R_0}
R(-a_j)$, and $t^{a_j}$ is short for the \textsc{Laurent} monomial
$t_1^{m_1} t_2^{m_2} \cdots t_n^{m_n}$ for $a_j = (m_1,\, m_2,\,
\cdots,\, m_n)$.

In view of\eqref{eq:iso}, applying the Additivity Theorem for
Characteristic Filtrations \cite[Corollary~2, p.~107]{Quillen-K} to
the filtration~\eqref{eq:filtration} shows that $\Phi$~has a
section~$\Xi$ given by the map whose components are induced by the
functors $P \mapsto T(P)_{a_j}$, $0 \leq j \leq m$. It is easy to see
that $\Xi$~is also right inverse to~$\Phi$; indeed, given $P \in
\pp(R_0)$ we see that $T \circ G^j (P) = P \tensor_{R_0} R(-a_j)
\tensor_R R_0$ is concentrated in degree~$a_j$, and hence canonically
isomorphic to~$P$ as an ungraded module.

Since $v$~was chosen in the interior of~$C$ we have $\bigcup_{k \in
  \bN} -kv+C = \bR^n$ so that $\ppgr(R)$ is the filtered union of the
$\ppgr_k(R)$. Since $K$\nbd-theory commutes with filtered unions
\cite[p.~104]{Quillen-K} the Theorem follows by considering the limit
process $k \to \infty$.

\medbreak

{\bf Examples.} The Theorem covers the case of coordinate rings of
affine toric schemes having a fixed point under the toric action, that
is, the case of any monoid ring of the form $R = S[C \cap \bZ^n]$
where $C$~is a full-dimensional pointed polyhedral rational cone
in~$\bR^n$. We may allow non-commutative rings~$S$ here as well as
non-rational cones~$C$, loosening the relation to toric
geometry. Specific examples include the polynomial ring $R_1 = S[X,Y]$
with the obvious $\bZ^2$-grading $\deg X = (1,0)$ and $\deg Y =
(0,1)$, and the ring $R_2 = S[U,V,W]/(UW-V^2)$ where $\deg U = (1,0)$,
$\deg(V) = (1,1)$ and $\deg(W) = (1,2)$; this can be realised using
for $C \subset \bR^2$ the cone generated by the vectors $(1,0)$
and~$(1,2)$. For a ring that is not a finitely generated
$S$\nbd-algebra consider the monoid ring $R_3 = S[C \cap \bZ^2]$ where
$C$~is the cone in~$\bR^2$ generated by the vectors $(1,0)$ and $(1,
\sqrt2)$. The Theorem then asserts that there are isomorphisms
\[K_i^\gr (R_j) \iso K_i(S) \tensor \bZ[t_1,\, t_1^{-1},\, t_2,\,
t_2^{-1}] \ , \qquad j=1,\, 2,\, 3 \ .\]

\medbreak

{\bf Acknowledgements.} The author would like to thank Judith Millar
and Roozbeh Hazrat for helpful discussions of graded $K$\nbd-theory.

\medskip\rightline{(03.11.2010, updated 11.01.2011)}

\providecommand{\bysame}{\leavevmode\hbox to3em{\hrulefill}\thinspace}
\providecommand{\MR}{\relax\ifhmode\unskip\space\fi MR }
\providecommand{\MRhref}[2]{%
  \href{http://www.ams.org/mathscinet-getitem?mr=#1}{#2}
}
\providecommand{\href}[2]{#2}

\end{document}